\documentstyle[12pt]{article}
\newtheorem{twierdz}{Theorem }[section]
\newtheorem{stwierdz}[twierdz]{Proposition }
\newtheorem{lemat}{Lemma }[section]
\newtheorem{definicja}{Definition }[section]
\newtheorem{przyklad}{Example}

\def\ot{\otimes}
\def\ra{\longmapsto}

\def\da{\downarrow}
\def\ca{{\cal A}}

\def\ce{{\cal E}}

\def\ide{id_{\ce \ot \ce}}

\def\1n{^{(1)}}
\def\2n{^{(2)}}
\def\3n{^{(3)}}

\def\t{\tilde}
\def\ba{\begin{array}}
\def\ea{\end{array}}
\def\be{\begin{equation}}
\def\ee{\end{equation}}
\def\bdm{\begin{displaymath}}
\def\edm{\end{displaymath}}
\def\beqnar{\begin{eqnarray}}
\def\eeqnar{\end{eqnarray}}
\def\beqna*{\begin{eqnarray*}}
\def\eeqna*{\end{eqnarray*}}

\renewcommand\thepage{\Large\bf May 1995 IFT UWr 890/\arabic{page}}
\begin{document}
\setcounter{page}{95}

\title{ON DEFORMATIONS OF COMMUTATION RELATION ALGEBRAS}
\author{Robert Ra{\l}owski\\
 \sc Institute of Theoretical Physics,
 \sc University of Wroc{\l}aw,\\
 \sc Pl. Maxa Borna 9,
 \sc 50-204  Wroc{\l}aw,\\
 \sc Poland\\
 \sc E-mail rober@ift.uni.wroc.pl}
\date{}
\maketitle
\vspace{2.5cm}
\begin{abstract}
This paper is on $C$ - symmetric creation
and annihilation operators, which are constructed on Wick's algebras
which
fulfil consistency conditions. The essential assumption is that every
algebraic action must be constant on equivalence classes. All
consistency
conditions follow from the above assumption. In this way we obtain
well
defined quotient algebras with some additional relations.
\end{abstract}
\vspace{3cm}
PACS numbers: 03.65.Ca, 03.65.Fd .

\thispagestyle{myheadings}

\newpage
\pagestyle{plain}
\renewcommand\thepage{\arabic{page}}
\setcounter{page}{1}

\section{Introduction.}

We would like to present the construction of deformed algebras of
commutation relation. An example of a deformed commutation relation
is an interpolation
of bosonic and fermionic statistics, e.g.
\cite{Boz1,Boz2,Green,Fiv,LPo}.
 The problem of additional relations between pairs of annihilation
operators and
between creation operators was considerated for example in
\cite{Jorg,Mar}.
The construction of annihilation operators is based on non commutative
differential calculus \cite{Bor1,Bor2,wz} as a result of quantum
groups, see \cite{Wor,Wor2,Pusz,Conn}.

In section {\bf \ref{definicje}} we applied the contraction notion
(evaluation \cite{wr1,wr2} ) defined in the algebraic way, see
\cite{oz}
on the tensor algebra of linear space E by taking generalised twist
$C$
between dual spaces of $E$ and $E$. This contraction satisfies $C$ -
deformed
Leibniz rule. In this section we introduced partial representation of
the creation and annihilation operators defined on algebra $TE$. In
this
representation we have commutation relations between creators and
annihilators.
In this way we have obtained Wick's algebra in which any sequences
consisting of
creation and annihilation operators can be arranged in Wick's way, see
\cite{Jorg}.
If we want to obtain additional relation, we have to divide the
algebra $TE$ by
${\cal N}$ graded ideal $J\subset TE$, generated by twist $B\in
End(E\ot E)$
and construct the subspace $J_2^*\subset E^*\ot E^*$, generated by
twist
$\t{B}\in End(E^*\ot E^*)$.
After this we have constructed representation of the annihilation and
creation
on the algebra $\ca=TE/J$, which is projected from the representation
defined
on the algebra $TE$. If we want to make the above construction from
the twists
$B,C$ they have to satisfy certain relations, from which it appears
that partial representation must
preserve the ideal $J$. In this way we have obtained non commutative
differential
calculus on the quantum plane see \cite{wz}. The relation between
annihilators
and the same relations between creators are necessary condition for
introduction
of inner product such that $<d^+x,y>=<x,dy>\;\; x,y\in \ca$. That is
why we have
introduced $\pi^*$ - invariant property for the contraction.
In the section {\bf \ref{relacje}} we present deformed commutation
relations:
$$
[a_i,a_j^+]_C=\delta_{i,j}\;\;[a^i,a^j]_B=0\;\;[a_i,a_j]_{\tilde{B}}=0
$$
as a consequence of definitions introduced in section {\bf
\ref{definicje}}.

To obtain representation on the $\ca$, the operators defined on
algebra $TE$ have to preserve the ideal contained in $TE$, which
leads us to
the consistency conditions between tensors $B,\tilde{B},C$.
The conditions which satisfy assumptions of the theorem
\ref{konsystencja2}
are enough to construct the representation defined on algebra $\ca$.
In this
case tensors $B,C$ satisfy generalised braided symmetry which are in
\cite{Bor1,Bor2}, where there are similar consistency conditions.
The contraction has $\pi^*$ - invariant property when tensors $B,C$
satisfy
assumptions of the theorem {\bf \ref{konsystencja3}}.
In section {\bf \ref{przyklady}} we present example which fulfil the
assumptions
of the theorems in section {\bf \ref{konsystencja}}. Of course bosons
and fermions satisfy
these assumption, however we has shown a q-deformed example.
\newpage

\section{Definitions and notions.}
\label{definicje}
Let us assume that $E$ is a  linear space of finite dimension over
field ${\cal C}$
with basis $\{f_i : i \in I\}$ and $E^*$ is the dual space to $E$
with
$\{e_j : j \in I\}$ dual basis.
Let $E^{0} := {\cal C}$ and $E^{\ot n} = E \ot E^{\ot (n-1)}$ be the
n-th tensorian
product of the space $E$.

\begin{em}
\begin{definicja}\label{Psi1}
Let $C :E^{*} \ot E \ra E \ot E^{*}$ be a map defined in the
following way:
\beqna*
 C (e_{i} \ot f_{j}) = \sum_{k,l} c_{i,j,k,l}\; f_{k} \ot e_{l}
 \; : c_{i,j,k,l} \in {\cal C }
\eeqna*
\end{definicja}
\end{em}

Now we define:
$C^{(i)}_n\in Hom(E^{\ot i-1}\ot (E^*\ot E)\ot E^{\ot n-i}),
E^{\ot i-1}\ot (E\ot E^*)\ot E^{\ot n-i})$
by
$$ C^{(i)}_n=1_{i-1}\ot C\ot 1_{n-i}\;\;for\;\;i\in\{1\ldots n\}$$
and
$$C^{(0)}=1=id\in End(E^*\ot E^{\ot n})$$
where $1_k=id_{E^{\ot k}}$. Let $ev: E^*\ot E\ra {\cal C}$ such that:
$$ev(e_i\ot f_j)=\delta_{i,j}$$ and extended by linearity to $E^*\ot
E$.

\begin{definicja}\label{ewaluacja}
The $C$ twisted contraction $ct_n$ with respect to given
elementary twist $C$ is a mapping:
$$ct_n:E^*\ot E^{\ot n}\ra E^{\ot n-1}$$
such that:
$$ct_n=\sum_{k=1}^n (1_{k-1}\ot ev\ot 1_{n-k})(C_n^{(k-1)}\ldots
C_n^{(0)})$$
\end{definicja}
Let us observe that $ct_1=ev$.

Let $A$ be any mapping :
$$E^*\ot E^{\ot n}\ni x\ot y \ra A_n(x\ot y)\in E^{\ot n-1}$$
and $W_i\in \{E,E^*\}$ for $i\in \{1\ldots k-1\}$ and $k\in {\cal
N}$, then we can define
$$A_{n,W_1\ldots W_{k-1}}^{(k)}\in Hom(W_1\ot\ldots\ot W_{k-1}\ot
E^*\ot E^{\ot n},
W_1\ot\ldots\ot W_{k-1}\ot E^*\ot E^{\ot n-1})$$
such that:
$$A_{n,W_1\ldots W_{k-1}}^{(k)}=1_{W_1\ldots W_k}\ot A_n
\;\;\;and\;\;\;k\in {\cal N}\;,$$
where $1_{W_1\ldots W_k}$ is the identity on the $W_1\ot\ldots\ot
W_{k-1}$.
In the natural way we can introduce an extension of $A_n^{(k)}$:
\beqna*
\tilde{A}_{n,W_1\ldots W_{k-1},W_1\ldots W_m}^{(k,m)}
&=&1_{W_1\ldots W_{k-1}}\ot A_n\ot 1_{W_1\ldots W_m}
\eeqna*

{}From now on we will indentify
$\tilde{A}_{n,W_1\ldots W_{k-1},W_1\ldots W_m}^{(k,m)}$ with
$A_{n,W_1\ldots W_{k-1}}^{(k)}$ becouse:
$$\tilde{A}_{n,W_1\ldots W_{k-1},W_1\ldots W_m}^{(k,m)}
=A_{n,W_1\ldots W_{k-1}}^{(k)}\ot id_{W_1\ot\ldots\ot W_m}$$
If $W_1=\ldots =W_{k-1}=E$, then we will accept the following
notation:
$A_n^{(k)}=A_{n,W_1\ldots W_{k-1}}^{(k)}$, if $W_1=E^*$ and
$W_2=\ldots =W_{k-1}=E$, then $A_{n,E^*}^{(k)}=A_{n,W_1\ldots W_{k-
1}}^{(k)}$.

Is easy to see that:
$$ct_n^{(1)}=\sum_{i=1}^n ct_1^{(i)}C^{(i-1)}\ldots C^{(1)}$$
where $ct_1^{(i)}=1_{i-1}\ot ct_1\ot 1_{n-i}$.

{}From the definition of contraction we obtain simple properties.
\begin{lemat}\label{l1}
The contraction has the following properties:
\begin{enumerate}
 \item contraction is a linear on $E^*\ot E^{\ot n}$
 \item $C$ - Leibniz rule:
$$ct_n^{(k)} = ct^{(k)}_{1} + ct^{(k+1)}_{n-1} \circ C^{(k)}
\;\;\;\;k\in {\cal N}$$
 \item $$ct_n^{(k)}=ct_1^{(k)}+\sum_{i=2}^n ct_1^{(i+k-1)}C^{(k+i-
2)}\ldots C^{(k)}$$
 \item $$ct_n^{(1)}=ct_m^{(1)}+ct_{n-m}^{(m+1)}C^{(m)}\ldots C^{(1)}$$
\end{enumerate}
\end{lemat}
{\sc Proof:} The point 1 is obvious, it remains to prove the points
2, 3.
Firstly we provide the proof this for $k=1$.
If $n=2$ then $C$ - Leibniz rule is simply satisfied. Then for $n>2$,
we have:
\beqna*
ct_{n+1}^{(1)}(e_i\ot f_j\ot y)&=&ct_1^1(e_i\ot f_j\ot
y)+\sum_{s=2}^{n+1}
ct_1^{(s)}\;C^{(s-1)}\ldots C^{(1)}(e_i\ot f_j\ot y)
\eeqna*
\beqna*
\hspace{1cm}&=&ct_1^1(e_i\ot f_j\ot y)+\sum_{k,l}\sum_{s=2}^{n+1}
c_{i,j,k,l}\;
ct_1^{(s)}\;C^{(s-1)}\ldots C^{(2)}(f_k\ot e_l\ot y)\\
&=&ct_1^1(e_i\ot f_j\ot y)+\sum_{k,l} c_{i,j,k,l}\;f_k\ot
\sum_{s=1}^{n}ct_1^{(s)}C^{(s-1)}\ldots C^{(1)}(e_l\ot y)\\
&=&ct_1^1(e_i\ot f_j\ot y)+\sum_{k,l} c_{i,j,k,l}\;f_k\ot
ct_n^{(1)}(e_l\ot y)\\
&=&ct_1^1(e_i\ot f_j\ot y)+\sum_{k,l} c_{i,j,k,l}\; ct_n^{(2)}(f_k\ot
e_l\ot y)\\
&=&ct_1^1(e_i\ot f_j\ot y)+ct_n^{(2)}C^{(1)}(e_i\ot f_j\ot y)
\eeqna*
and for arbitrary $k\in {\cal N}$, $w_i\in W_i$, where
$i\in \{1\ldots k-1\}$, $x\in E^*,\;\;y\in E^{\ot n}$ :
\beqna*
ct_n^{(k)}(w_1\ot\ldots \ot w_{k-1}\ot x\ot y)&=&
(w_1\ot\ldots \ot w_{k-1}\ot ct_n^{(1)}(x\ot y))
\eeqna*
\beqna*
\hspace{2cm}&=&(w_1\ot\ldots \ot w_{k-1}
\ot (ct_1^{(1)}+ct_{n-1}^{(2)}\;C^{(1)})(x\ot y))\\
&=&(ct_1^{(k)}+ct_{n-1}^{(k+1)}\;C^{(k)})
(w_1\ot\ldots \ot w_{k-1}\ot x\ot y)
\eeqna*

Now we prove the point 3 as a simple consequence of definition
\ref{ewaluacja} :
\beqna*
ct_n^{(k)}(w_1\ot\ldots \ot w_{k-1}\ot x\ot y)&=&
(w_1\ot\ldots \ot w_{k-1}\ot ct_n^{(1)}(x\ot y))
\eeqna*
\beqna*
\hspace{2cm}&=&(w_1\ot\ldots \ot w_{k-1}\ot \sum_{i=1}^n
ct_1^{(i)}C^{(i-1)}\ldots C^{(1)}(x\ot y))\\
&=&(ct_1^{(k)}+\sum_{i=2}^n ct_1^{(k+i-1)}\;C^{(k+i-2)}\ldots
C^{(k)})(w_1\ot\ldots \ot w_{k-1}\ot x\ot y)
\eeqna*
The point 4 follows from the $3^{rd}$ point.
This completes the proof\hfill $\Box$

It is interesting to note that $C$ - Leibniz rule and $A_1$ gives the
corresponding
$A_n\in Hom(E^*\ot E^{\ot n},E^{\ot n-1})$ in unique way.
\begin{lemat}\label{l2}
If $A_n\in Hom(E^*\ot E^{\ot n}, E^{\ot n-1})$ such that:
\begin{enumerate}
 \item $$A_1x\ot y=X(y)$$
 \item $\exists C\in Hom(E^*\ot E,E\ot E^*)$ and
    $$A_n^{(1)}=A_1^{(1)}+A_{n-1}^{(2)}\;C^{(1)}\;\;\;\;for\;\;n\in
\{ 2,3,\ldots\}$$
\end{enumerate}
then
$$A_n^{(1)}=\sum_{i=1}^n A_1^{(i)}C^{(i-1)}\ldots C^{(1)}$$
\end{lemat}
{\sc Proof:} For $n=2$ we have:
$$A_2^{(1)}=A_1^{(1)}+A_1^{(2)}C^{(1)}$$
Let us suppose that for $n=m$ the lemma is true, then:
\beqna*
A_{m+1}^{(1)}(e_i\ot f_j\ot y)&=&(A_1^{(1)}+A_m^{(2)}C^{(1)})(e_i\ot
f_j\ot y)\\
&=&A_1^{(1)}(e_i\ot f_j\ot y)+\sum_{k,l}c_{i,j,k,l}\;A_m^{(2)}(f_k\ot
e_l\ot y)\\
&=&A_1^{(1)}(e_i\ot f_j\ot y)+\sum_{k,l}c_{i,j,k,l}\;f_k\ot
A_m^{(1)}(e_l\ot y)\\
&=&A_1^{(1)}(e_i\ot f_j\ot y)+\sum_{k,l}c_{i,j,k,l}\;f_k\ot
\sum_{s=1}^m A_1^{(s)}\;C^{(s-1)}\ldots C^{(1)}(e_l\ot y)\\
&=&A_1^{(1)}(e_i\ot f_j\ot y)+\sum_{k,l}c_{i,j,k,l}\;
\sum_{s=1}^m A_1^{(s+1)}\;C^{(s)}\ldots C^{(2)}(f_k\ot e_l\ot y)\\
&=&A_1^{(1)}(e_i\ot f_j\ot y)+
\sum_{s=1}^m A_1^{(s+1)}\;C^{(s)}\ldots C^{(2)}(C(e_i\ot f_j)\ot y)\\
&=&A_1^{(1)}(e_i\ot f_j\ot y)+
\sum_{s=1}^m A_1^{(s+1)}\;C^{(s)}\ldots C^{(2)}C^{(1)}(e_i\ot f_j\ot
y)\\
&=&\sum_{s=1}^{m+1} A_1^{(s)}\;C^{(s-1)}\ldots C^{(1)}(e_i\ot f_j\ot
y)
\eeqna*
The proof is completed\hfill $\Box$

\begin{em}
\begin{definicja}\label{represent}
The operators $a_{n,i}$ and $a_{n,i}^{+}$ for $n\in {\cal N}$we call
respectively partial
annihilation and creation $C$ - operators, when we have:
\beqna*
a_{n,i}(y)&:=&ct_n^{(1)}(e_{i}\ot y)\;\;\;\;y\in E^{\ot n}\\
a_{n,j}^{+}(y)&:=&f_{j}\ot y\;\;\;\;y\in E^{\ot n}
\eeqna*
\end{definicja}
\end{em}

Let $B\in Hom(E \ot E,E \ot E)$. We can define the subspace in $E \ot
E$
by the following relation:
$$J_{2} = \{ x \in E^{\ot 2} |\; \exists_{z\in E^{\ot 2}}\;
x= z - B(z) \}=Im(1-B)$$
And more general:
$$J_{n}=E\ot J_{n-1}+J_{n-1}\ot E \; \; \mbox{ for } n=3,4,...$$

So we can divide $E^{\ot n}$ by $J_{n}$ obtaining linear subspace
$A^{\odot n}:=E^{\ot n}/J_{n}$

Let
\bdm J=\bigoplus_{n\in {\cal N}} J_{n}\edm
This is a minimal ideal in the tensor algebra $TE$ containing $J_{2}$.
In this case we can simply introduce quotient algebra
$\ca=TE/J$, which is a ${\cal N}$ graded
space:$A=\bigoplus_{n>=0}A^{\odot n}$\\
We have canonical homomorfism $\pi_{n}$, which projects space $E^{\ot
n}$ onto
$A^{\odot n}$. For the second order it gives:
$$E^{\ot 2} \ni x \ot y \ra \pi_{2}(x \ot y) = [x \ot y] \in A^{\odot
2}$$
and $\pi_2$ satisfies the relation:
\be
\pi_{2} = \pi_{2} \circ B
\label{diag1a}
\ee

Every element of $A^{\odot n}$ is in image of the projector $\pi_{n}$
:
$$for\;\;\;every\;\;\;y\in A^{\odot n} \;\;there\;\; exists\;\; y_{1}
\in E^{\ot n}
\;\;such\;\; that\;\; y=\pi_{n}(y_{1})$$

So, we can introduce the canonical projector $\pi$ of the space $TE$
onto $\ca$ space:
\be
\pi : TE \ra \ca
\label{canon1}
\ee

Let $\tilde{B}\in End(E^{*}\ot E^{*})$, then we can define $J_2^*$
in the following way:  $$J_2^*=Im (1-\t{B})$$

Now we introduce $B$ - commutator:
\begin{em}
$$[a_{n+1,i}^{+},a_{n,j}^{+}]_B:=a_{n+1,i}^{+} \circ a_{n,j}^{+} -
\sum_{k,l}b_{i,j,k,l}
a_{n+1,k}^{+} \circ a_{n,l}^{+}$$
\end{em}
and also $\tilde{B}$ and $C$ commutators:
\begin{em}
$$[a_{n-1,i},a_{n,j}]_{\tilde{B}}:=a_{n-1,i} \circ a_{n,j} -
\sum_{k,l}\t{b}_{i,j,k,l}
a_{n-1,k} \circ a_{n,l}$$
$$[a_{n+1,i},a_{n,j}^+]_{C}:=a_{n+1,i} \circ a_{n,j}^+ -
\sum_{k,l}c_{i,j,k,l}
a_{n-1,k}^+ \circ a_{n,l}$$
where $b_{i,j,k,l}\;\t{b}_{i,j,k,l}\;c_{i,j,k,l}$ are matrix elements
of
$B$, $\t{B}$, $C$ respectively.
\end{em}

To complete these definitions, we would like to introduce a
definition of $\pi^*$ - invariant mapping on the whole space $TE$,
which will be play important role in Proposition \ref{relat2} below.

\begin{em}
\begin{definicja}\label{niezm_a}
Let $w\in Hom(E^{*} \ot TE, TE)$ lowering the level of tensor by 1,
then we will say that $w$ is $\pi^{*}$ - inwariant on $TE$ when:
$$w_{n-1}^{(1)}\circ w_{n,E^*}^{(2)}\;(J_2^*\ot E^{\ot n})\subset
J_{n-2}\;\;\;for\;\;n=(2,3,\ldots )$$
where $w_n\in (E^*\ot E^{\ot n},E^{\ot n-1})$ are the corresponding
components of $w$.
\end{definicja}
\end{em}

\begin{em}
\begin{definicja}\label{pirepr}
Let $d_{n,i}:A^{\odot n}\ra A^{\odot (n-1)}\;\;
and\;\;d_{n,i}:A^{\odot n}\ra A^{\odot (n+1)}$
for $n\in {\cal N}$, then $(d_{n,i},d_{n,i}^{+})$ is $\pi^*$ -
invariant
$C$ - representation of the annihilation and creation operators, when
exists
the $C$ - partial representation $a_{n,i}\;a_{n,i}^+$, which
satisfies the
following conditions:
\begin{enumerate}
 \item \label{diag2}
   The operators $d_{n,i}$ fulfil the following diagram 
$$\ba{ccc}
&a_{n,i}&\\
E^{\ot n}&\ra&E^{\ot n-1}\\[0.2 cm]
\pi_n\;\da& &\da\;\pi_{n-1}\\[0.2 cm]
A^{\odot n}&\ra&A^{\odot n-1}\\
&d_{n,i}&
\ea$$
 \item \label{diag3}
  The operators $d_{n,i}^+$ fulfil the following diagram 
$$\ba{ccc}
&a_{n,i}^+&\\
E^{\ot n}&\ra&E^{\ot n+1}\\[0.2 cm]
\pi_n\;\da& &\da\;\pi_{n-1}\\[0.2 cm]
A^{\odot n}&\ra&A^{\odot n+1}\\
&d_{n,i}^+&
\ea$$
 \item The operator $a_i$ given by ( $a_{n,i}(x)=ct_n^{(1)}(e_i\ot
x)$ ), is $\pi^{*}$ - invariant
  operator see definition \ref{niezm_a}
\end{enumerate}
\end{definicja}
\end{em}
If we have operators $d_{n,i},\;d_{n,i}^+$ defined on $A^{\odot n}$
then we can
construct operators $d_i,\;d_j^+$ defined on $\ca$ algebra in the
following
way: $$d_i=\bigoplus_{n\in {\cal N}}d_{n,i}\hspace{3cm}
d_i^+=\bigoplus_{n\in {\cal N}}d_{n,i}^+$$

\section{On $C$ and $B$ - relations.}
\label{relacje}
In this section we would like to calculate $C$ commutation relation
in the representation introduced in section {\bf \ref{definicje}}.

Firstly we show the theorem given by J\"{o}rgensen, Schmith and
Werner \cite{Jorg} :

\begin{em}
\begin{twierdz}[J\"{o}rgensen, Schmith, Werner]
 Let $(a_{i},a_{j}^{+})$ be $C$ - partial representation
then we have the following relation:
\bdm
[a_{i},a_{j}^{+}]_{C} = \delta_{i,j} \; id
\label{relat1}
\edm
\end{twierdz}
\end{em}

{\sc Proof: } Let $y \in E^{\ot n}$ for every $n \in N$, then by
application
of the $C$ - Leibniz rule we have:
\beqna*
a_{n+1,i} \circ a_{n,j}^{+}(y) & = & a_{n+1,i}(f_{j} \ot y)\\
& = & ct_{n+1}^{(1)}(e_{i} \ot f_{j} \ot y)\\
& = & ct_{1}^{(1)}(e_{i} \ot f_{j} \ot y)+
ct_{n}^{(2)} \circ C^{(1)}(e_{i} \ot f_{j} \ot y )\\
& = & \delta_{i,j}\; y + \sum_{k,l}c_{i,j,k,l}\;(f_{k} \ot
(ct_{n}(e_{l} \ot y)))\\
& = & [\delta_{i,j} +\sum_{k,l}c_{i,j,k,l}\;(a_{n-1,k}^{+} \circ
a_{n,l})](y)\\
\label{dow1}
\eeqna*
So the proof is completed $\hfill \Box$

In quotient space $\ca$ we have additional relations generated by
ideal $J \subset TE$. So we have the following two propositions:

\begin{em}
\begin{stwierdz}\label{relat2}
 If $d_{i},\; \; d_{j}$ are $C$ - annihilation
operators which are $\pi^*$ - invariant (see Definition \ref{pirepr})
then
\bdm
[d_{i},d_{j}]_{\tilde{B}} = 0.
\edm
\end{stwierdz}
\end{em}

{\sc Proof: } Let $y_{1} \in A^{\odot n}$ and $y_{1} = \pi_{n}(y)$
for $y \in E^{\ot n}$.
Using definition \ref{pirepr} point \ref{diag2}
and $\pi^*$ - invariantness see definition \ref{niezm_a}, we have
simple calculation:
\beqna*
d_{n-1,i} \circ d_{n,j}(y_{1}) & = & \pi_{n-2} \circ a_{n-1,i} \circ
a_{n,j}(y)\\
&=&\pi_{n-2}ct_{n-1}^{(1)}ct_{n}^{(2)}(e_{i} \ot e_{j} \ot y)\\
&=&\pi_{n-2}ct_{n-1}^{(1)}ct_{n}^{(2)}(e_{i} \ot e_{j}-
\tilde{B}(e_{i}\ot e_{j}))\ot y\\
&&+ \pi_{n-2}ct_{n-1}^{(1)}ct_{n}^{(2)}\tilde{B}(e_{i}\ot e_{j})\ot
y\\
&=&\pi_{n-2}z+\sum_{k,l}\t{b}_{i,j,k,l}\;\pi_{n-2}ct_{n-
1}^{(1)}ct_{n}^{(2)}(e_k\ot e_l\ot y)\\
&=&\sum_{k,l}\t{b}_{i,j,k,l}\;d_{n-1,k}\circ d_{n,l}(y_1)\\
\label{dowod2}
\eeqna*
Where $z\in J_{n-2}$. That completed the proof $\hfill \Box$

\begin{em}
\begin{stwierdz}\label{relat3}
 If $d_{i}^{+},\; \; d_{j}^{+}$ are $C$ - creation operators
which are $\pi^*$ - invariant (see Def \ref{pirepr}) then we have the
following
relation:
\bdm
[d_{i}^{+},d_{j}^{+}]_{B^{*}} = 0.
\edm
\end{stwierdz}
\end{em}

{\sc Proof: } Let $y_{1} \in A^{\odot n}$ and $y_{1} = \pi_{n}(y)$
for $y \in E^{\ot n}$.
Using equation (\ref{diag1a}) and definition \ref{pirepr} point
\ref{diag3},
we calculate:
\beqna*
d_{n+1,i}^{+} \circ d_{n,j}^{+}(y_{1}) & = & \pi_{n+2} \circ
a_{n+1,i}^{+} \circ a_{n,j}^{+}(y) \\
& = & \pi_{n+2}(f_{i} \ot f_{j} \ot y)\\
& = & \pi_{n+2}(f_{i} \ot f_{j}-B(f_{i} \ot f_{j}))\ot y\\
&&+\pi_{n+2}B(f_{i} \ot f_{j}) \ot y\\
&=&\pi_{n+2}z+\sum_{k,l}b_{i,j,k,l}\;\pi_{n+2}(f_{k}\ot f_{l}\ot y)\\
&=&\sum_{k,l}b_{i,j,k,l}\;d_{n+1,k}\circ d_{n,l}(y_1)
\label{dowod3}
\eeqna*
Where $z\in J_{n+2}$. That completes the proof $\hfill \Box$

\section{Consistency conditions.}
\label{konsystencja}
To obtain well defined algebra with additional multiplication
relations,
we have to guarant the constant action on the layers in algebra $TE$.
For example every contraction defined in above section must be well
defined
on the layers in $TE$, hence we have the following proposition:

\begin{em}
\begin{stwierdz}\label{consist}
The necessary condition for the constant action on the layers of
$a_i$ operator
is:
$$(id - B)(id + \tilde{C}) = 0$$
where $\t{C}$ has the matrix elements given by
$\tilde{c}_{i,j,k,l}=c_{j,l,i,k}$.
\end{stwierdz}
\end{em}

{\sc Proof:} By taking elements of the basis of the space $E$, $f_{j}
\; f_{k}$ and
the dual space $E^{*}$ and its basis $e_{i}$ and using the equation
(\ref{diag1a}) and
$C$ - Leibniz rule see lemma \ref{l1} we have the simple calculation:
\beqna*
\bar{0}&=&ct_2(e_i\ot (1-B)(f_j\ot f_k))\\
&=&\sum_{l,m}(\delta_{j,l}\delta_{k,m}-b_{j,k,l,m})ct_2(e_i\ot f_l\ot
f_m)\\
&=&\sum_{l,m}(\delta_{j,l}\delta_{k,m}-
b_{j,k,l,m})(ct_1^{(1)}+ct_1^{(2)}C^{(1)})
(e_i\ot f_l\ot f_m)\\
&=&\sum_{l,m}(\delta_{j,l}\delta_{k,m}-
b_{j,k,l,m})(\delta_{i,l}f_m+\sum_{r,s}
c_{i,l,r,s}ct_1^{(2)}\;(f_r\ot e_s\ot f_m))\\
&=&\sum_{l,m}(\delta_{j,l}\delta_{k,m}-b_{j,k,l,m})
\sum_r(\delta_{i,l}\delta_{m,r}\;f_r+c_{i,l,r,m}\;f_r)\\
&=&\sum_r\sum_{l,m}(\delta_{j,l}\; \delta_{k,m} - b_{j,k,l,m})\;
(\delta_{i,l}\;\delta_{m,r}+c_{i,l,r,m})\;f_r
\eeqna*
Then we have:
$$\sum_{l,m}(\delta_{j,l}\; \delta_{k,m} - b_{j,k,l,m})\;
(\delta_{i,l}\;\delta_{m,r}+c_{i,l,r,m})=0$$
$\hfill \Box$

Above result is the only necessary condition to obtain well defined
operator $d_i,d_i^{+}$ on
quotient algebra $\ca$. This is well know as linear Wess Zumino
condition see
\cite{wz}.
We also know the following:
\begin{em}
\begin{stwierdz}
The operators $d_i$ exist, (see Definition \ref{pirepr})
iff for every $n\in {\cal N}$ the following relation hold:
$$\pi_{n}(x) = \pi_{n}(y) \Rightarrow \pi_{n-1} \circ a_{n,i}(x) =
\pi_{n-1} \circ a_{n,i}(y)$$
\end{stwierdz}
\end{em}

Now we will present the proof of the theorem given by A. Borowiec and
V. Kharchenko
see \cite{Bor1}:
\begin{em}
\begin{twierdz}[Borowiec-Kharchenko]\label{konsystencja2}
If the tensors $B\in End(E\ot E)$, and $C\in Hom(E^*\ot E,E\ot E^*)$
satisfy the following conditions:
 \begin{enumerate}
  \item $(1-B)(1+\tilde{C})=0$ where
$\tilde{c}_{i,j,k,l}=c_{j,l,i,k}$\label{pierwszy}
  \item there exists $A\in Hom(E^*\ot E^{\ot 2},E^{\ot 2}\ot E^*)$
such that: $\;C^2C^1B^2-B^1C^2C^1=(1-B^1)A$ \label{drugi}
\label{trzeci}
 \end{enumerate}
then $a_r(J)\subset J$ where $J\subset TE$ is the minimal ideal
generated by the relation\\
$<x-Bx=0>$ for $x\in E\ot E$.
\end{twierdz}
\end{em}

{\sc Proof:} At first we have to prove for $J_{2}$ and $J_{3}$.
Condition for the ideal $J_2$ is described in the proposition
\ref{consist}.
We have to prove that for every $r\in I\;\;\;a_{3,r}(J_3)\subset J_2$.
Then
\beqna*
LHS&=&a_{3,r}((1-B)y_1\ot x_1 + x_2\ot (1-B)y_2) \in J_2
\eeqna*
for every $y_1,y_2\in E\ot E$, and $x_1,x_2\in E$, so we have:
\beqna*
\hspace{-0.2cm}
LHS&=&(ct_1^{(1)}+ct_1^{(2)}C^1+ct_1^{(3)}C^{(2)}C^{(1)})((1-
B^{(2)})e_r\ot y_1\ot x_1+(1-B^{(3)})e_r\ot x_2\ot y_2)
\eeqna*
\beqna*
\hspace{2cm}&=&(ct_1^{(1)}+ct_1^{(2)}C^{(1)})((1-B^{(2)})e_r\ot
y_1\ot x_1+(1-B^{(3)})e_r\ot x_2\ot y_2)\\
&+&ct_1^{(3)}C^{(2)}C^{(1)}((1-B^{(2)})e_r\ot y_1\ot x_1+(1-
B^{(3)})e_r\ot x_2\ot y_2)
\end{eqnarray*}
{}From the assumption \ref{pierwszy} and \ref{drugi} we have:
\beqna*
(ct_1^{(1)}+ct_1^{(2)}C^{(1)})(1-B^{(2)})(e_r\ot y_1\ot x_1)=0\in E\\
ct_1^{(3)}C^{(2)}C^{(1)}B^{(2)}=Bct_1^{(3)}C^{(2)}C^{(1)}+(1-
B)ct_1^{(3)}A
\eeqna*
Then we have:
\beqna*
LHS&=&(ct_1^{(1)}+ct_1^{(2)}C^{(1)})(1-B^{(3)})e_r\ot x_2\ot y_2+(1-
B)ct_1^3C^2C^1e_r\ot y_1\ot x_1\\
\eeqna*
\beqna*
\hspace{2cm}&+&(1-B)ct_1^{(3)}A(e_r\ot y_1\ot
x_1)+ct_1^{(3)}C^{(2)}C^{(1)}(1-B^{(3)})e_r\ot x_2\ot y_2
\eeqna*
The second and the third terms belong to $J_2$, so we have:
\beqna*
LHS&=&(ct_1^{(1)}+ct_1^{(2)}C^{(1)}+ct_1^{(3)}C^{(2)}C^{(1)})(1-
B^{(3)})e_r\ot x_2\ot y_2+z\\
&=&ct_3^{(1)}(1-B^{(3)})e_r\ot x_2\ot y_2+z
\eeqna*
where $z\in J_2$, but
\beqna*
ct_3^{(1)}(1-B^{(3)})e_r\ot x_2\ot
y_2&=&(ct_1^{(1)}+ct_2^{(2)}C^{(1)})(1-B^{(3)})e_r\ot x_2\ot y_2
\eeqna*
\beqna*
\hspace{2cm}&=&e_r(x_2)\;(1-B)y_2+\sum_s\alpha_s ct_2^{(2)}C^{(1)}(1-
B^{(3)})e_r\ot f_s\ot y_2\\
&=&e_r(x_2)\;(1-B)y_2+\sum_s\sum_{p,t}\alpha_s\;c_{r,s,p,t}f_p\ot
ct_2^{(1)}
(1-B^{(2)})e_t\ot y_2\\
&=&e_r(x_2)\;(1-B)y_2+\sum_p f_p\ot \vec{0}
\eeqna*
then $LHS$ is in the $J_2$.

If for $n=(2,3)$ the theorem is true, then we will show this for
every $n>3$ by
induction. So, let's suppose that $ct_{n-1}(E^*\ot J_{n-1})\subset
J_{n-2}$.
Every element of $J_n$ we can write as:
$$y=\sum_{k=1}^{n-1}(1-B^{(k)})y_k\;\;\;\; \forall_{k\in
\{1..k\}}\;\;y_k\in E^{\ot n}$$
then for all $x\in E^* $ we have:
\beqna*
ct_n^{(1)}\; x\ot y&=&\sum_{k=1}^{n-1}ct_n^{(1)}\;(x\ot (1-
B^{(k)})y_k)\\
&=&\sum_{k=1}^{n-1}ct_n^{(1)}\;(1-B^{(k+1)})\;x\ot y_k\\
&=&ct_n^{(1)}(1-B^{(2)})x\ot y_1+ct_n^{(1)}(1-B^{(3)})\;x\ot y_2\\
&&+\sum_{k=3}^{n-1}ct_n^{(1)}\;(1-B^{(k+1)})\;x\ot y_k
\eeqna*
but
\beqna*
\sum_{k=3}^{n-1}ct_n^{(1)}\;(1-B^{(k+1)})\;x\ot y_k&=&
\sum_{k=3}^{n-1}(ct_1^{(1)}+
ct_{n-1}^{(2)}C^{(1)})\;(1-B^{(k+1)})\;x\ot y_k\\
&=&\sum_{k=3}^{n-1}(1-B^{(k)})\tilde{y'}_k\\
&&+\sum_{k=3}^{n-1}ct_{n-1}^{(2)}(1-B^{(k)})
\sum_{i,j=1}^{dim(E)}\tilde{y}_i\ot \tilde{x}_j\ot \tilde{y}_{k,i,j}
\eeqna*
of course the first and the second terms belong to $J_{n-1}$
from the assumption.
\beqna*
ct_n^{(1)}(1-B^{(2)})&=&(ct_1^{(1)}+ct_{n-1}^{(2)}C^{(1)})\;(1-
B^{(2)})\\
&=&(ct_1^{(1)}+ct_1^{(2)}C^{(1)}+ct_{n-2}^{(3)}C^{(2)}C^{(1)})\;(1-
B^{(2)})\\
&=&(ct_2^{(1)}+ct_{n-2}^{(3)}C^{(2)}C^{(1)})\;(1-B^{(2)})\\
&=&ct_2^{(1)}(1-B^{(2)})+ct_{n-2}^{(3)}(1-B^{(1)})C^{(2)}C^{(1)}
+ct_{n-2}^{(3)}(1-B^{(1)})A^{(1)}\\
&=&{\vec{0}}+(1-B^1)ct_{n-2}^{(3)}\;C^{(2)}C^{(1)}
+(1-B^{(1)})ct_{n-2}^{(3)}A^{(1)}
\eeqna*
\beqna*
ct_n^{(1)}(1-
B^{(3)})&=&(ct_1^{(1)}+ct_1^{(2)}C^{(1)}+ct_1^{(3)}C^{(2)}C^{(1)}+ct_{}\cr
&=&ct_3^{(1)}(1-B^{(3)})+ct_{n-3}^{(4)}(1-
B^{(2)})C^{(3)}C^{(2)}C^{(1)}\\
&&+(1-B^{(2)})ct_{n-3}^{(4)}A^{(2)}C^{(1)}
\eeqna*
Then
\beqna*
ct_n^{(1)}(1-B^{(2)})x\ot y_1&&\in J_{n-1}\\
ct_n^{(1)}(1-B^{(3)})x\ot y_2&&\in J_{n-1}
\eeqna*
That completed the proof \hfill $\Box$

\begin{em}
\begin{twierdz}\label{konsystencja3}
If matrixes $\t{B}\in End(E^*\ot E^*), C\in Hom(E^*\ot E,E\ot E^*)$
satisfy the following conditions:
 \begin{enumerate}
  \item fulfil YBE:
$\t{B}^{(2)}C^{(1)}C^{(2)}=C^{(1)}C^{(2)}\t{B}^{(1)}$\label{first}
  \item $[ct_1^{(1)}ct_{1,E^*}^{(k)}C^{(k-1)}\ldots C^{(2)}+ct_1^{(k-
2)}C^{(k-3)}\ldots C^{(1)}ct_{1,E*}^{(2)}][1-\t{B}^{(1)}]=0$
  for $k\in \{3..n\}\\ and\;\;\;\; n\in {\cal N}$\label{second}
 \end{enumerate}
then $ct_{n-1}^{(1)}\;ct_{n,E^*}^{(2)}(J_2^*\ot E^{\ot n})\;\subset
J_{n-2}\subset E^{\ot (n-2)}$
\end{twierdz}
\end{em}
{\sc Proof:} We have to prove this by induction. Let us suppose that
for certain
$n-1\in {\cal N}$ thesis is satisfied, then we have to prove that
assertion is
satisfied for $n$. Then
\beqna*
ct_{n-1}^{(1)}ct_{n,E^*}^{(2)}(1-\t{B}^{(1)})&=&(ct_1^{(1)}+ct_{n-
2}^{(2)}C^{(1)})(ct_{1,E^*}^{(2)}+ct_{n-1,E^*}^{(3)}C^{(2)})(1-
\t{B}^{(1)})
\eeqna*
\beqna*
\hspace{1.5cm}&=&(ct_1^{(1)}ct_{n,E^*}^{(2)}+ct_{n-
2}^{(2)}C^{(1)}ct_{1,E^*}^{(2)})(1-\t{B}^{(1)})+ct_{n-
2}^{(2)}C^{(1)}ct_{n-1,E^*}^{(3)}C^{(2)}(1-\t{B}^{(1)})\\
&=&(ct_1^{(1)}ct_{n,E^*}^{(2)}+ct_{n-
2}^{(2)}C^{(1)}ct_{1,E^*}^{(2)})(1-\t{B}^{(1)})+ct_{n-2}^{(2)}ct_{n-
1,E^*}^{(3)}C^{(1)}C^{(2)}(1-\t{B}^{(1)})\\
&=&(ct_1^{(1)}ct_{n,E^*}^{(2)}+ct_{n-
2}^{(2)}C^{(1)}ct_{1,E^*}^{(2)})(1-\t{B}^{(1)})+ct_{n-2}^{(2)}ct_{n-
1,E^*}^{(3)}(1-\t{B}^{(2)})C^{(1)}C^{(2)}
\eeqna*
The last term from assertion belongs to the $J_{n-2}$. Then
\beqna*
ct_{n-1}^{(1)}ct_{n,E^*}^{(2)}(1-
\t{B}^{(1)})&=&(ct_1^{(1)}ct_{n,E^*}^{(2)}+ct_{n-
2}^{(2)}C^{(1)}ct_1^{(2)})(1-\t{B}^{(1)})+z
\eeqna*
\beqna*
\hspace{2cm}&=&(ct_1^{(1)}ct_{n,E^*}^{(2)}+ct_{n-
2}^{(2)}C^{(1)}ct_{1,E^*}^{(2)})(1-\t{B}^{(1)})\\
&&+(ct_1^{(1)}ct_{1,E^*}^{(2)}-ct_1^{(1)}ct_{1,E^*}^{(2)})(1-
\t{B}^{(1)})+z\\
&=&(ct_1^{(1)}ct_{n,E^*}^{(2)}+ct_{n-1}^{(1)}ct_{1,E^*}^{(2)}-
ct_1^{(1)}ct_{1,E^*}^{(2)})(1-\t{B}^{(1)})+z\\
&=&ct_1^{(1)}(ct_{1,E^*}^{(2)}+ct_1^{(3)}C^{(2)}+\ldots
+ct_{1,E^*}^{(n+1)}C^{(n)}\ldots C^{(2)})\\
&&+(ct_1^{(1)}+ct_1^{(2)}C^{(1)}+\ldots +ct_1^{(n-1)}C^{(n-2)}\ldots
C^{(1)})ct_{1,E^*}^{(2)})(1-\t{B}^{(1)})\\
&&-ct_1^{(1)}ct_{1,E^*}^{(2)}(1-\t{B}^{(1)})+z
\eeqna*
where $z\in J_{n-2}$. The assertion we have from the
$\ref{second}^{nd}$
assumption.

To finish the proof we have to calculate for $n=2$ only:
\beqna*
ct_1^{(1)}ct_{2,E^*}^{(2)}(1-\t{B}^{(1)})&=&
ct_1^{(1)}(ct_{1,E^*}^{(2)}+ct_{1,E^*}^{(3)}\;C^{(2)})(1-
\t{B}^{(1)})\\
&=&(ct_1^{(1)}ct_{1,E^*}^{(2)}+ct_1^{(1)}ct_{1,E^*}^{(3)}\;C^{(2)})(1-
\t{B}^{(1)})\;\;=\;\;0
\eeqna*
\hfill $\Box$

\section{Examples.}\label{przyklady}
In this section we would like to show a few examples, which fulfil
consistency
conditions formulated in the previous section. We start with two
simple examples:
\begin{przyklad}[Bosons]
\bdm
b_{i,j,k,l}=\t{b}_{i,j,k,l}=c_{i,j,k,l}=\delta_{i,l}\;\;\delta_{j,k}
\edm
\end{przyklad}
\begin{przyklad}[Fermions]
\bdm
b_{i,j,k,l}=\t{b}_{i,j,k,l}=c_{i,j,k,l}= -
\;\delta_{i,l}\;\;\delta_{j,k}
\edm
\end{przyklad}

In general case linear condition is satisfied when
$Im(\ide -B)\subset Ker(\ide+C)$. This relation is not satisfied when:
$$B(f_i\ot f_j)=q^{i-j}f_j\ot f_i\hspace{1.5cm}\t{B}(e_i\ot e_j)=q^{i-
j}e_j\ot e_i
\hspace{1.5cm}C(e_i\ot f_j)=q^{i-j}f_j\ot e_i$$
This relation satisfies the $\ref{drugi}^{nd}$ assumption
of the theorem \ref{konsystencja2}. For $q\in \{-1,1\} $ all the
assumptions
of theorems \ref{konsystencja2} and \ref{konsystencja3} are
satisfied. Then
we have the following:
\begin{przyklad}[Mixed Bosons - Fermions]
$$b_{i,j,k,l}=\t{b}_{i,j,k,l}=c_{i,j,k,l}=(-1)^{i-
j}\;\;\delta_{i,l}\;\;\delta_{j,k}$$
\end{przyklad}
and more general:
\begin{przyklad}[q-Deformed algebra]
$$b_{i,j,k,l}=\t{b}_{i,j,k,l}=q^{j-
i}\;\delta_{i,l}\;\delta_{j,k}\hspace{2cm}
c_{i,j,k,l}=q^{i-j}\;\delta_{i,l}\;\delta_{j,k}\hspace{1cm} for\;\;\;
q\in {\cal R}$$ as particular case of more general relations, which
was studied in \cite{Mar}.
\end{przyklad}
\vspace{1cm}

\noindent {\bf Acknowledgements}\\
We would like to thank to R. Gielerak, A. Borowiec, M. Bo{\.z}ejko and
W. Marcinek for discussions and critical remarks.



\begin{thebibliography}{}
\bibitem{Boz1} M. Bo{\.z}ejko, R. Speicher, An example of a
Generalized Brownian
Motion, Comm. Math. Phys. {\bf 137} (1991), 519-531.
\bibitem{Boz2} M. Bo{\.z}ejko and R. Speicher. Completely positive
maps on
Coxeter groups, deformed commutation relations and operator spaces.
Preprint Heidelberg, Summer 1993.
\bibitem{Green} O.W. Greenberg, Particles with small violations of
Fermi or
 Bose statistics, Phys. Rev. {\bf D 43} (1991), 4111-4120.
\bibitem{Fiv} D. Fivel, Interpolation between Fermi and Bose
statistics using
generalised commutators, Phys. Rev. Lett. {\bf 65} (1990), 3361-3364.
\bibitem{LPo} R. Lenczewski and K. Podg{\'o}rski, A q-analog of the
Quantum
Central Limit Theorem for $SU_q( 2)$, J. Math. Phys. {\bf 33} (1992),
(2768-2778).
\bibitem{Jorg} P.E.T. J\"{o}rgensen, L.M. Schmith, R.F. Werner,
Positive
 representation of general commutation relations allowing Wick
ordering.
 Preprint 1993.
\bibitem{Mar} W. Marcinek, On unital braidings and quantization.
Preprint
ITP UWr No 847 (1993), to be published in  Rep. Math. Phys.
\bibitem{Bor1} Andrzej Borowiec and Vladislav K. Kharchenko,
Free calculi and coordinate algebras.
\bibitem{Bor2} Andrzej Borowiec, Vladislav K. Kharchenko, Zbigniew
Oziewicz,
On free differentials on associative algebras.
In  "Non-Associative algebras and its applications", edited by Santoz
Gonzales,
Kluwer Academic Publishers, Dordrecht 1994 [Mathematics and its
Applications,
vol 303], pp 46-53.
\bibitem{wz} J. Wess and B. Zumino, Covariant differential calculus
on the
quantum hyperplane. CERN preprint 5697/90.
\bibitem{Wor} S. Woronowicz and W. Pusz, Twisted second quantization,
Rep. Math. Phys. {\bf 27} (1989), 231-257.
\bibitem{Wor2} S. Woronowicz, Differential calculus on compact matrix
pseudogroups (quantum groups), Commun. Math. Phys. 122 (1989), 125-
170.
\bibitem{Pusz} W. Pusz, Rep. Math. Phys. {\bf 27} (1989), 394
\bibitem{Conn} A. Connes, Non-commutative differential geometry,
Publ.Math.
IHES {\bf 62} (1985), 257-360.
\bibitem{wr1} W. Marcinek R. Ra{\l}owski, Particle operators from
braided
geometry. Proceeding of the XXX Winter School of Theoretical Physics
in Karpacz, Poland, February 15 - 26, (1994).
\bibitem{wr2} W. Marcinek R. Ra{\l}owski, On Wick algebra's with
braid relations.
Preprint ITP UWr No 876 (1994), to be published in J. Math. Phys.
\bibitem{oz} Z. Oziewicz, Cz. Sitarczyk, Parallel treatment of
Riemannian
and symplectic Clifford algebra's. A. Micali et al., Clifford
algebra's
and their applications in mathematical physics, 83-95.
 1992 Kluwer Academic Publishers. Printed in the Netherlands.
\end{thebibliography}
\end{document}